\documentclass[12pt, draft]{amsart}
\usepackage[cp850]{inputenc}
\usepackage{amssymb}
\usepackage[mathscr]{eucal}
\usepackage{amscd}
\usepackage{amsmath}
\usepackage[all]{xy}

\newcommand{\R}{\mathbb R}
\newcommand{\N}{\mathbb N}

\newcommand{\C}{\mathbb C}

\newcommand{\KP}{{\sf K}\hspace{-1pt}{\sf P}}

\newcommand{\Rn}{\mathfrak R^{(n)}}

\theoremstyle{plain}
\newtheorem{theorem}{Theorem}[section]
\newtheorem{prop}[theorem]{Proposition}
\newtheorem{lemma}[theorem]{Lemma}
\newtheorem{cor}[theorem]{Corollary}
\newtheorem{definition}[theorem]{Definition}

\newtheorem{defn}[theorem]{Definition}

\newcommand{\pin}[2]{\langle #1, #2 \rangle}
\newcommand{\di}{\mathrm{dim}}
\newcommand{\Sp}{\mathrm{Sp}}
\newcommand{\GL}{\mathrm{GL}}

\newcommand{\adef}{\begin{defn}}
\newcommand{\zdef}{\end{defn}}

\begin{document}

\title{Symplectic forms on Banach spaces}

\author[J.M.F. Castillo]{Jes\'us M. F. Castillo}
\address{Universidad de Extremadura, Instituto de Matem\'aticas Imuex, E-06011 Badajoz, Spain.} \email{castillo@unex.es}

\author[W. Cuellar]{Wilson Cuellar }
\address{Departamento de Matem\'atica, Instituto de Matem\'atica e
Estat\'istica, Universidade de S\~ao Paulo, rua do Mat\~ao 1010,
05508-090 S\~ao Paulo SP, Brazil}
\email{cuellar@ime.usp.br}

\author[M. Gonz\'alez]{Manuel Gonz\'alez}
\address{Departamento de Matem\'aticas,
Universidad de Cantabria, E-39071 Santander, Spain.}
\email{manuel.gonzalez@unican.es}

\author[R. Pino]{Ra\'ul Pino}
\address{Departamento de Matem\'aticas,
Universidad de Extremadura, E-06011 Badajoz, Spain.}
\email{rpino@unex.es \hspace*{2mm} raul.pino@rai.usc.es}

\thanks{2020 Mathematics Subject Classification. Primary: 46B20, 46B10; Secondary  46M18, 46B70\\
Keywords: Symplectic Banach space; symplectic operator;Rochberg spaces; Kalton-Peck space; Hilbert space}

\thanks{The research of the first and third authors was supported in part by MINCIN project PID2019-103961GB. The research of the first author was supported in part by Junta de Extremadura project IB20038. The research of the second author was supported by FAPESP grants  (2016/25574-8), (2018/18593-1) and (2019/23669-0). The research of the fourth author was partially supported by project FEDER-UCA18-108415 funded by 2014-2020 ERDF Operational Programme and by the Department of Economy, Knowledge, Business and University of the Regional Government of Andalucia.}

\maketitle
\thispagestyle{empty}

\begin{abstract} We extend and generalize the result of Kalton and Swanson ($Z_2$ is a symplectic Banach space with no Lagrangian subspace) by showing that all higher order Rochgberg spaces $\mathfrak R^{(n)}$ are symplectic Banach spaces with no Lagrangian subspaces. The nontrivial symplectic structure on even spaces is the one induced by the natural duality; while the nontrivial symplectic structure on odd spaces requires perturbation with a complex structure. We will also study symplectic structures on general Banach spaces and, motivated by the unexpected appearance of complex structures, we introduce and study almost symplectic structures. \end{abstract}

\begin{flushright} This paper is dedicated\\ to the memory of Yuliia Zdanovska,\\ brilliant mathematician promise\\ killed in Jarkhov by Putin's war.\\ Long may live her Teach for Ukraine project\end{flushright}

\section{Introduction}

A  real Banach space $X$ is said to be \emph{symplectic} if there is a  continuous alternating bilinear map $\omega: X\times X \to \R$ such that the induced map $L_\omega \, : \, X \to X^*$ given by $L_\omega(x)(y)=\omega(x,y)$ is an isomorphism onto.  A symplectic Banach space is necessarily isomorphic to its dual and reflexive (see Lemma \ref{symref}). During the decade of 1970's several authors drew the attention to the importance of the study of symplectic forms on Banach spaces and, more broadly, on  Banach manifolds.  For instance, in the proof of Weinstein  \cite{wein} of an infinite dimensional version of the classical Darboux theorem for symplectic geometry, or in the Hamiltonian formulation of infinite dimensional mechanics due to  Chernoff and Marsden \cite{chermar}. See also Swanson \cite{Swanson, Swan1} for various results about symplectic structures on Banach spaces.\medskip

A motivation for this work has been  the negative solution given by Kalton and Swanson \cite{KalSwan} to the question raised by Weinstein \cite{wein} of whether every infinite dimensional symplectic Banach space is trivial. A symplectic Banach space $(X, \omega)$  is said \emph{trivial}  if there does exist a reflexive Banach space $Y$ and an isomorphism $T: X \to Y\oplus Y^*$ such that $\omega(x,y)= \Omega_Y(Tx,Ty)$ for every $x,y\in X$, where
$$\Omega_Y[ (z,z^*), (w,w^*)] = w^*(z) - z^*(w).$$
In this case,  $T^{-1} (Y\times \{0\})$ is, according to Definition \ref{lagrangian}, a \emph{Lagrangian} subspace of $(X, \omega)$.
Kalton and Swanson show that the celebrated Kalton-Peck space $Z_2$ (see \cite{Kalton_Peck}) is a symplectic space with no Lagrangian subspaces. In this paper we will consider  the sequence of higher order Rochberg spaces $\mathfrak R^{(n)}$ \cite{rochberg} obtained from the scale of $\ell_p$ spaces, which can be considered as generalizations of both $\ell_2$ and $Z_2$ since $\mathfrak R^{(1)}= \ell_2$ and $\mathfrak R^{(2)}= Z_2$. We will show that all these spaces $\mathfrak R^{(n)}$ are symplectic and contain no (infinite dimensional) Lagrangian subspaces; in other words, they admit a nontrivial symplectic structure. A remarkable point is that while the nontrivial symplectic structure on even spaces  $\mathfrak R^{(2n)}$ is the one induced by the natural duality; the nontrivial symplectic structure on the odd spaces  $\mathfrak  R^{(2n+1)}$ requires to modify the natural duality structure with a complex structure.\medskip

A second motivation for this work is to clarify the connection between symplectic and complex structures on Banach spaces. Recall that given a real Banach space  $X$, a linear and bounded operator $J: X \to X$ is called a \emph{complex structure} on $X$ if $J^2=-Id$. In this case, the operator $J$ induces a $\C$-linear structure on $X$ in the form  $i x = J(x)$. In the Hilbert space setting there is a correspondence between symplectic structures and complex structures since Weinstein \cite[Prop. 5.1]{wein} proved that  every symplectic structure on a Hilbert  space $\mathcal H$  has the form $\omega (x,y) = \pin{J(x)}{y}$, for some complex structure $J$ and an equivalent inner product on $\mathcal H$.  The argument of Weinstein actually shows that every symplectic structure on a real Hilbert space is trivial.
This correspondence is no longer valid in general (say, non-reflexive spaces may admit complex structures), but still some properties of complex structures can be studied in the context of symplectic structures. The section \ref{almost} is devoted to study perturbations of symplectic structures by strictly singular operators and extensions of symplectic structures on hyperplanes following the techniques of V. Ferenczi \cite{F} and Ferenczi and E. Galego \cite{FGEven} about  complex structures. We also prove an analogous result for symplectic structures to those of \cite{ccfm} for complex structures: no symplectic structure on $\ell_2$ can be extended to a bilinear form on a hyperplane $H$ of $Z_2$ containing it.

\section{Background}

 \begin{definition}
 Given a real  Banach space $X$, a \emph{linear symplectic  form} on $X$ is a  bilinear map $\omega \, : \, X\times X\to  \R$  satisfying the following

\begin{enumerate}
\item $\omega $ is continuous: there exists $K>0$ such that $|\omega(x,y)|\leq K \|x\| \|y\|$ for every $x,y\in X$.
\item  $\omega$ is alternating: $\omega(x,y)=-\omega(y,x)$ for all $x,y\in X$.
\item The induced map $L_\omega \, : \,X\to X^*$ given by $L_\omega (x)(y)=\omega(x,y)$ is an isomorphism of $X$ onto $X^*$.
\end{enumerate}
In this case, the pair $(X, \omega)$ of a Banach space with a symplectic form is called a symplectic Banach space. We say that $X$ is \emph{symplectic} if there exists a linear symplectic  form $\omega$ on $X$.
\end{definition}

The following result is from the pioneering work of Kalton and Swanson \cite{KalSwan}
\begin{lemma}\label{symref} A continuous alternating bilinear map $\omega$ on a real Banach space $X$ is symplectic if and only if $X$ is reflexive and $L_\omega:X\to X^*$ is an isomorphism into.\end{lemma} We include the proof for the sake of completeness:

\begin{proof}
If $(X, \omega)$ is symplectic, then $L_\omega^* \, : \, X^{**} \to X^*$ is also an isomorphism onto and if $x\in X\subseteq X^{**}$ then $L_\omega ^*(x)= -L_\omega(x)$.  So $X=X^{**}$ in the canonical embedding. Assume that $X$ is reflexive and that $L_\omega:X\to X^*$ is an isomorphism into.  Suppose that there exists $f\in X^* \setminus \mathrm{Im} (L_\omega)$. By the Hanhn-Banach theorem in combination with reflexivity there exists $y\in X$ such that $L_\omega (x) (y)=0$ for every $x\in X$ and  $f(y)\neq 0$. Then $\omega(x,y)=0$  for every $ x \in X$ and so $L_\omega(y) =0$ which contradicts that $L_\omega$ is injective.\end{proof}

Thus, a real  Banach space $X$ is symplectic if and only if $X$ is reflexive and there exists $\alpha:X\to X^*$  isomorphism into such that $\alpha^*=-\alpha$ (where $\alpha^*: X^{**} \to X^*$ is the adjoint of $\alpha$ with the canonical identification). In this case $\omega(x,y) =\alpha(x)(y)$ is a symplectic form on $X$. These results justify the following definition.

\begin{definition} Let $X$ be a real reflexive Banach space. An isomorphism $\alpha:X\to X^*$  is said to be symplectic if $\alpha^*=-\alpha$. \end{definition}

The basic examples of symplectic Banach spaces known so far are:
\begin{itemize}
\item Finite dimensional symplectic spaces are even dimensional.
\item If $X$ is a  real reflexive  Banach space and $$ \omega ( (e,e^*), (f,f^*))= f^*(e)-e^*(f)
$$ then $(X\oplus X^*, \omega)$ is a symplectic space.
\item Infinite dimensional Hilbert spaces are symplectic.
\item The Kalton-Peck twisted Hilbert spaces $\ell_{2}(\phi)$ introduced in \cite{Kalton_Peck} are symplectic.
\end{itemize}

Only symplectic structures on Hilbert spaces admit a simple description.

\begin{lemma} \label{hilbertsymplectic} Let $(\mathcal H, \pin{\,}{\,})$ be a real Hilbert space and $\omega$ be a continuous bilinear form on $\mathcal H$. Then $\omega$ is a symplectic form on $\mathcal H$ if and only if  there exists an isomorphism $\alpha: \mathcal H \to \mathcal H$ with Hilbert-adjoint $\alpha^{*} = -\alpha$ such that $\omega(x,y)=\pin{\alpha(x)}{y}$ for all $x,y\in \mathcal H$.
\end{lemma}

We pass to define the equality notion for symplectic forms.

\begin{definition}\label{defsymp} Two symplectic spaces $(X_1, \omega_1)$ and $(X_2, \omega_2)$ are equivalent if there is an onto isomorphism $T: X_1 \to X_2$ such that $\omega_2(Tx,Ty)= \omega_1(x,y)$. If $\omega_1, \omega_2$ are two symplectic structures on a space $X$ we will write
$\omega_1\sim \omega_2$ to denote they are equivalent.\end{definition}

Hence the symplectic spaces $(X_1, \omega_1)$ and $(X_2, \omega_2)$  are equivalent if and only if there exists an isomorphism $T: X_1 \to X_2$  such that $L_{\omega_1} = T^* L_{{\omega}_2}T$, where $L_{{\omega}_i}$ is the respective isomorphism of $X_i$ onto $X_i^*$, $i=0,1$.
$$\xymatrixrowsep{1cm}
\xymatrixcolsep{2cm}
\xymatrix{X_1 \ar[r]^T \ar[d]_{L_{\omega_1}} & X_2 \ar[d]^{L_{\omega_2}} \\
X_1^*  &  X_2^*\ar[l]^{T^*}}$$

The following notions will also be useful:

\begin{definition} \label{lagrangian} Let $(X,\omega)$ be a symplectic space and let  $F$ be a closed subspace  of $X$.

\begin{itemize}
\item The symplectic orthogonal  (or symplectic annihilator) of $F$ is the linear subspace $F^\omega= \{ x\in X\, : \, \omega(x,y)=0 \, \, \, \text{for all} \, \, \, y\in F\}$.
\item $F$ is symplectic if $(F, \omega|_{F\times F})$ is symplectic.
\item $F$ is isotropic if $\omega(x,y)=0$ for every $x,y\in F$, i.e., $F\subseteq F^{\omega}$.
\item  $F$ is Lagrangian if it is isotropic and possesses an isotropic complement.
\item The symplectic structure $(X,\omega)$ is trivial if there is a Lagrangian subspace $F$ such that $X\simeq F\oplus G$ and $(X,\omega)$ is equivalent to $(F\oplus G, \omega|_F\times \omega|_G )$.
\end{itemize}
\end{definition}

Observe that for a closed subspace $F=(F^{\omega})^{\omega}$. 
By Zorn's lemma every symplectic structure admits a maximal isotropic subspace that must be necessarily closed.

\begin{lemma}  \label{symsubspace}Let $(X,\omega)$ be a real  symplectic space. A closed subspace   $F$  of $X$ is symplectic if and only if $X=F\oplus F^\omega$.
\end{lemma}
\begin{proof}  Denote by $L_F=L_{\omega|_{F\times F}}$, then $L_F$ is injective if and only if $F\cap F^{\omega}=\{0\}$. Suppose that $F$ is symplectic. For every $x\in X$  there exists $f\in F$ such that $L_F(f)=L_\omega(x)_{|F}$. It follows that $x-f\in F^\omega$ and hence $X=F\oplus F^\omega$.  Conversely, assume that $X=F\oplus F^\omega$. Given $\phi\in F^*$, let $x\in X$ such that $\phi= L_\omega(x)_{|F}$. If we write $x=f+g$ with $f\in F$ and $g\in F^\omega$, then $L_F(f)=\phi$.  We conclude that $L_F$ is an isomorphism and hence $F$ is symplectic.
\end{proof}

From here it immediately follows:

\begin{cor}\label{symio}\label{codimension} Let $(X,\omega)$ be a real symplectic space. A closed subspace   $F$  of $X$ is symplectic if and only if $F^\omega$ is symplectic. If a closed finite codimensional subspace $F$ of $X$ is symplectic then the dimension of $X/F$ is even.
\end{cor}

As an immediate consequence, no real symplectic structure extends from an hyperplane to the whole space. Moreover

\begin{lemma}\label{subsym} Let $(X,\omega)$ be a real  symplectic space and let $H$ be a closed hyperplane of $X$. Then  $H^\omega \subseteq H$ and $\di (H^\omega)=1$. Moreover, there exists a closed symplectic subspace $H'\subseteq H$ with $\di(H/H') =1$.
\end{lemma}
\begin{proof} Let $g\in X^*$ with $H= \ker g$. Since $X$ is symplectic, there exists $x_0\neq 0$ in $X$ such that $L_\omega(x_0)=g$. For every $h\in H$, we have $\omega(x_0,h)=L_\omega(x_0)(h)=g(h)=0$. Hence $x_0\in H^\omega$. On the other hand, $g(x_0)=\omega(x_0,x_0)=0$, and then $x_0\in H$. Now for every $x\in H^\omega$, we have $\ker g \subseteq \ker L_\omega(x)$. It follows that  $L_\omega(x) = \lambda L_\omega(x_0)$  for some constant $\lambda\in \R$. We conclude that $x=\lambda x_0$ and then $\di (H^\omega)=1$. For the moreover part, let $X= H\oplus \mathrm{span} \{x_1\}$. The argument above implies that there exists a closed hyperplane $H_1$ of  $X$ containing $x_1$ such that $H=\{x_1\}^\omega$.  Consider $H'= H \cap H_1$. It follows that $X=  H' \oplus \mathrm{span} \{x_0,x_1\}$ and since $\omega (x_0,x_1) \neq 0$ we have that $ \{x_0,x_1\}^\omega =H'$. We conclude that $H'$ is symplectic by using Lemma  \ref{symsubspace}. \end{proof}

\section{About the symplectic structure of $Z_2$}

The only known ``nontrivial" symplectic spaces are the twisted Hilbert spaces $\ell_2(\phi)$; among them, the most remarkable and notorious is the Kalton-Peck space $Z_2$. Let us review the main features of $Z_2$ to then advance in the study of its symplectic structure.
The space $Z_2$ is generated by the Kalton-Peck map $\KP x = x \log x$ in the form $Z_2= \ell_2\oplus_{\KP }\ell_2 = \{(\omega, x)\in \ell_\infty \times \ell_2: \omega - \KP x\in \ell_2\}$ endowed with the quasinorm (but equivalent to a norm) $\|(\omega, x)\| = \|\omega - \KP x\|_2 + \|x\|_2$. The space $Z_2$ is a twisted Hilbert space in the sense that there is an exact sequence
$$\xymatrix{0\ar[r]&\ell_2\ar[r]^{\jmath}& Z_2 \ar[r]^{Q}&\ell_2\ar[r]&0}$$ with inclusion $\jmath(y) = (y,0)$ and quotient map $Q(\omega, x)= x$. The space $Z_2$ is isometric to its dual (see below) but is not isomorphic to a Hilbert space; it admits a basis but not an unconditional basis. It however has an Unconditional Finite Dimensional Decomposition into the $2$-dimensional subspaces $X_n={\rm{span}}\{(e_n,0),(0,e_n)\}$. An operator $T\in \mathfrak L(Z_2)$ is either strictly singular or an isomorphism on a complemented copy of $Z_2$.  An operator $T\in \mathfrak L(Z_2)$ is strictly singular if and only if $T\jmath$ is strictly singular \cite{Kalton_Peck}.

The most remarkable fact for our purposes is that $Z_2$ is isometric to its dual. Indeed  \cite[Th. 5.1]{Kalton_Peck}$Z_2 = \ell_2 \oplus_{\KP} \ell_2$ and $Z_2^* = \ell_2\oplus_{-\KP}\ell_2$. The duality between $Z_2$ and $Z_2^*$  is given by \cite{kaltsym}
$$ \langle (x,y),(x',y') \rangle = \langle x,y' \rangle_{\ell_2} + \langle y,x' \rangle_{\ell_2}.$$
See \cite{CabeCast2021} for details. This yields that any of the maps $(x,y)\to (-x,y)$ or $(x,y)\to (x, -y)$ is an isometry $Z_2\to Z_2^*$. Accordingly,

\begin{prop}\label{P_01} $Z_2$ is a symplectic Banach space. Precisely, the bilinear antisymmetric form $\lhd \cdot, \cdot \rhd:Z_2\times Z_2\rightarrow \mathbb{R}$ given by
$$ \lhd (x,y),(x',y') \rhd =\langle x,y' \rangle_{\ell_2} - \langle y,x' \rangle_{\ell_2}.$$
is such that $D:Z_2\rightarrow Z_2^*$ given by $D(a)[b]=\lhd\, a,b \,\rhd$ is an isomorphism.
\end{prop}

\adef\label{eq_02} The symplectic adjoint $T^+:Z_2\rightarrow Z_2$ of an operator $T\in \mathfrak L(Z_2)$ is defined assigning to each $y\in Z_2$ the only vector $T^+y$ such that, for all $x\in Z_2$,
\begin{equation}
 \lhd\, T^+y,x \,\rhd = \lhd\, y,Tx \rhd.
\end{equation}\zdef

Indeed, $T^+$ exists since the map $x\to \lhd\, y,Tx \,\rhd$ defines a continuous functional on $Z_2$. By Proposition \ref{P_01}, there exists an unique $y'\in Z_2$ so that $\lhd\, y',x \,\rhd=\lhd\, y,Tx \,\rhd$ for all $x\in Z_2$. Observe that in the complex case the map $D:Z_2\rightarrow Z_2^*$ must be made antilinear  as in the case of duality between Hilbert spaces and one can define an involution $+:\mathcal{B}(Z_2)\rightarrow\mathcal{B}(Z_2)$ such that $T^+$ identifies with the usual dual map $T^*$. The map $T^+$ is bounded whenever $T$ is bounded since $\frac{1}{3}\|T\|\leq \|T^+\|\leq 3\|T\|$, which can be proved using that $\|D\|\leq 3$ and $\|D^{-1}\|\leq 1$ (see \cite{CabeCast2021}). Moreover, there is a commutative diagram
\begin{equation}\label{diagrama_02}
\xymatrixrowsep{1cm}
\xymatrixcolsep{2cm}\xymatrix{Z_2 \ar[r]^D \ar[d]_{T^+} & Z_2^* \ar[d]^{T^*} \\
Z_2 \ar[r]^{D} &  Z_2^*}
\end{equation}
where $D:Z_2\rightarrow Z_2^*$ is the isomorphism given in Proposition \ref{P_01}. Indeed, given $x\in Z_2$, one has $DT^+x=\lhd\, T^+x, \cdot \,\rhd\in Z_2^*$ and the other way around gives $T^*Dx=T^*\big(\lhd\, x, \cdot\, \rhd\big)=\lhd\, x,T(\cdot) \,\rhd\in Z_2^*$, and both functionals coincide by (\ref{eq_02}). This duality was fully exploited by Kalton \cite{kaltsym} and Kalton and Swanson \cite{KalSwan}.

\subsection{Matrix representation for $T^+$}

Bounded operators $T:X\rightarrow X$ defined on reflexive Banach spaces with basis $(e_i)_{i\in\mathbb{N}}$ admit a matrix representation
$(a_{ij})$ in the sense that $T(e_i)=\sum_{j=1}^\infty a_{ij}e_j$ in terms of such basis. Indeed, the canonical duality between $X$ and $X^*$ given by $\langle e_i,e_j^*\rangle=\delta_{ij}$ yields
 $a_{ij}= \langle T(e_i),e_j^* \rangle$. Taking into account the identities
$\langle T^*(e_i^*),e_j^{**}\rangle=\langle T^*(e_i^*),e_j\rangle=T^*(e_i^*)(e_j)=e_i^*(Te_j)=\langle e_i^*,T(e_j)\rangle=\langle T(e_j),e_i^*\rangle$ it is then clear that the matrix representation of $T^*$ is just the transpose of that of $T$.

$Z_2$ is a superreflexive Banach space \cite{CabeCast2021} with a basis $(u_n)_{n\in\mathbb{N}}$ defined for each $n\in\mathbb{N}$ by $u_{2n-1}=(e_n,0)\quad\text{and}\quad u_{2n}=(0,e_n)$, where $(e_n)_n$ is the canonical basis of $\ell_2$ \cite[Th. 4.10]{Kalton_Peck}. The symplectic form $\lhd\,\cdot,\cdot\,\rhd$ above defines the matrix $(\lhd\,u_i, u_j\,\rhd)_{ij}$, namely
\begin{equation}\label{eq_04}
\begin{pmatrix}
\lhd\, u_1 , u_1 \,\rhd & \lhd\, u_1 , u_2 \,\rhd & \lhd\, u_1 , u_3 \,\rhd & \dots\\
\lhd\, u_2 , u_1 \,\rhd & \lhd\, u_2 , u_2 \,\rhd & \lhd\, u_2 , u_3 \,\rhd & \dots\\
\lhd\, u_3 , u_1 \,\rhd & \lhd\, u_3 , u_2 \,\rhd & \lhd\, u_3 , u_3 \,\rhd & \dots\\
\vdots & \vdots & \vdots & \ddots\\
\end{pmatrix}=\begin{pmatrix}
0 & 1 & 0 & 0 & \cdots\\
-1 & 0 & 1 & 0 & \cdots\\
0 & -1 & 0 & 1 & \cdots\\
0 & 0 & -1 & 0 & \cdots\\
\vdots & \vdots & \vdots & \vdots & \ddots\\
\end{pmatrix}.
\end{equation}
On the other hand, any operator $T:Z_2\rightarrow Z_2$ admits a matrix representation $(a_{ij})$ so that $T(u_i)=\sum_{j=1}^\infty a_{ij}u_j$.
Taking (\ref{eq_04}) into account we deduce that the matrix $(a_{ij})$ is
\begin{align*}
\begin{pmatrix}
\lhd\, T(u_1) , u_2 \,\rhd & -\lhd\, T(u_1) , u_1 \,\rhd & \lhd\, T(u_1) , u_4 \,\rhd & -\lhd\, T(u_1) , u_3 \,\rhd & \dots \\
\lhd\, T(u_2) , u_2 \,\rhd & -\lhd\, T(u_2) , u_1 \,\rhd & \lhd\, T(u_2) , u_4 \,\rhd & -\lhd\, T(u_2) , u_3 \,\rhd & \dots\\
\lhd\, T(u_3) , u_2 \,\rhd & -\lhd\, T(u_3) , u_1 \,\rhd & \lhd\, T(u_3) , u_4 \,\rhd & -\lhd\, T(u_3) , u_3 \,\rhd & \dots\\
\lhd\, T(u_4) , u_2 \,\rhd & -\lhd\, T(u_4) , u_1 \,\rhd & \lhd\, T(u_4) , u_4 \,\rhd & -\lhd\, T(u_4) , u_3 \,\rhd & \dots\\
\vdots & \vdots & \vdots & \vdots &\ddots
\end{pmatrix}.
\end{align*}

Then, $\lhd\, T(u_i),\cdot\,\rhd$ corresponds to the row $i$, and $\lhd\,\cdot,u_j\,\rhd$ corresponds to column $j+1$ or $j-1$ whether $j$ is, respectively, odd or even. It follows that
$$a_{ij}=(-1)^{j+1}\lhd\, T(u_i), u_{j+(-1)^{j+1}} \,\rhd$$
Now, if $b_{ij}=(-1)^{j+1}\lhd\, T^+(u_i), u_{j+(-1)^{j+1}} \,\rhd$ are the coefficients of the matrix representation of $T^+$, reasoning in the same manner we obtain
\begin{align*}
b_{ij}&=(-1)^{j+1}\lhd\, T^+(u_i), u_{j+(-1)^{j+1}} \,\rhd=(-1)^{j+1}(-1) \lhd\, T(u_{j+(-1)^{j+1}}), u_{i} \,\rhd\\
&=(-1)^j\lhd\, T(u_{j+(-1)^{j+1}}), u_{i} \,\rhd = (-1)^{j+i} a_{j+(-1)^{j+1}\,i+(-1)^{i+1}}
\end{align*}

Summing up  $T^+\equiv(b_{ij})=\big((-1)^{i+j}a_{j+(-1)^{j+1}\,i+(-1)^{i+1}}\big)$, i.e.,

$$ T=\begin{pmatrix}
a_{11} & a_{12} & a_{13} & a_{14} & \dots\\
a_{21} & a_{22} &  a_{23} & b_{24} & \dots\\
a_{31} & a_{32} & a_{33} & a_{34} & \dots\\
a_{41} & a_{42} & a_{43} & a_{44} & \dots\\
\vdots & \vdots & \vdots & \vdots & \ddots\\
\end{pmatrix}
\Longrightarrow T^+ = \begin{pmatrix}
a_{22} & -a_{12} & a_{42} & -a_{32} & \dots\\
-a_{21} & a_{11} & -a_{41} & a_{31} & \dots\\
a_{24} & -a_{14} & a_{44} & -a_{34} & \dots\\
-a_{23} & a_{13} & -a_{43} & a_{33} & \dots\\
\vdots & \vdots & \vdots & \vdots & \ddots\\
\end{pmatrix}\medskip.$$

The paper \cite{cgr} contains a study of operators on $Z_2$ with a different matrix representation: observe that every operator on $Z_2$ can be represented with a matrix $\begin{pmatrix}
\alpha & \beta\\
\delta & \gamma
\end{pmatrix}$ where $\alpha, \beta, \gamma, \delta$ are linear maps $\mathbb K^\N\to \mathbb K^\N$. It is easy to check now that if $T=\begin{pmatrix}
\alpha & \beta\\
\delta & \gamma
\end{pmatrix}$ is a bounded operator on $Z_2$ then $T^+= \begin{pmatrix}
\gamma^* & -\beta^*\\
-\delta^* & \alpha^*
\end{pmatrix}$.
The matrix $\begin{pmatrix}0& 1\\
0 & 0
\end{pmatrix}$ represents the bounded operator $Z_2 \stackrel{Q}\to \ell_2 \stackrel{\jmath}\to Z_2$ and $\begin{pmatrix}
0 & 1\\
0 & 0
\end{pmatrix}^+=\begin{pmatrix}
0& -1\\
0 & 0
\end{pmatrix}.$
Consequently, $(\jmath Q)^+ (\jmath Q) = (\jmath Q)^2= 0$, so the $C^*$-algebra identity $\|T^+T\|=\|T\| \,\|T^+\|$ fails and thus $Z_2$ can not be renormed in such a way that $(\mathfrak L(Z_2),+)$ becomes a $C^*$-algebra (even if one redefines the involution for the complex case). This is related to the classical Kawada-Kakutani-Mackey Theorem \cite{Kakutani_Mackey,Kawada} because the involution $+$ is (in his terms) \emph{not proper} (see \cite{Spain}).

\section{Symplectic transformations on $Z_2$}

\adef An operator $T:Z_2\rightarrow Z_2$ will be called a \emph{symplectic transformation} if it preserves the symplectic form, in the sense that \begin{equation}\label{eq_05}
\lhd\, T(x),T(y) \,\rhd=\lhd\, x,y \,\rhd,\quad\text{for all }x,y\in Z_2.
\end{equation}
\zdef

An operator $T$ is a symplectic transformation if and only if $T^+T=I$ (here $I$ is the identity): indeed, for all $x,y\in Z_2$ we have
$$\lhd\,T^+Tx,y\,\rhd-\lhd\,x,y\,\rhd=\lhd\,(T^+T-I)x,y\,\rhd=0,$$
and thus we deduce from Proposition \ref{P_01} that $(T^+T-I)x=0$ for all $x\in Z_2$. The other implication is clear. From this we obtain:

\begin{prop} Symplectic transformations in $Z_2$ have complemented range.
\end{prop}

Unbounded symplectic transformations on $Z_2$ are possible: just set the linear map $L(e_n, 0)= (e_n, 0)$ and $L(0, e_n)=(ne_n, e_n)$. Indeed, $L$ preserves $\lhd\, \cdot,\cdot \,\rhd$ by checking on the basis elements, and it is unbounded since $\| L(e_n, e_n)\|=\|(n+1)e_n - \KP(e_n)\|_2 +\|e_n\|_2= n+2$, for every $n\in \mathbb N$. Let us show some natural bounded examples.

\adef An operator $\eta:\ell_2\to \ell_2$ is said to be an \emph{operator on the scale} if there is $p>2$ such that both $\eta: \ell_p\to \ell_p$ and $\eta: \ell_{p^*}\to \ell_{p*}$ are bounded. It will be called an isometric operator on the scale if both $\eta: \ell_p\to \ell_p$ and $\eta: \ell_{p^*}\to \ell_{p*}$ are into isometries.\zdef

A result of Banach \cite{banach} establishes that $U\big((x_n)_n\big)=(\varepsilon_nx_{\sigma(n)}),$ where $\sigma:\mathbb{N}\rightarrow\mathbb{N}$ is a permutation and $|\varepsilon_n|=1$ for all $n\in\mathbb{N}$ are the only examples of surjective isometric operators on the scale. One of the forms of the Commutator Theorem, see \cite{tams,racsam} for details, is that if $\eta$ is an operator on the scale then $\tau_\eta = \begin{pmatrix}
\eta & 0\\
0& \eta
\end{pmatrix}$ defines an operator on $Z_2$; see also \cite{cgr}. Operators of the form $\begin{pmatrix}
\alpha & \beta\\
0& \gamma\end{pmatrix}$ will be called upper triangular operators.

\begin{prop}\label{opT} An upper triangular operator is a symplectic transformation if and only if it has the form $\begin{pmatrix}
\alpha & S\\
0 & \alpha
\end{pmatrix}$ with $\alpha \in  \mathfrak L(\ell_2)$ an isometric operator on the scale and $S\in \mathfrak L(\ell_2)$ such that $\alpha^* S$ is selfadjoint.
\end{prop}
\begin{proof} $
\begin{pmatrix}
\alpha & S\\
0& \alpha
\end{pmatrix}^+ \begin{pmatrix}
\alpha & S\\
0& \alpha
\end{pmatrix} =
\begin{pmatrix}
\alpha^* & - S^*\\
0& \alpha^*
\end{pmatrix}\begin{pmatrix}
\alpha & S\\
0 & \alpha
\end{pmatrix}= \begin{pmatrix}
\alpha^*\alpha & \alpha^*S-S^*\alpha\\
0 & \alpha^*\alpha
\end{pmatrix}=\begin{pmatrix}
1 & 0\\
0 & 1
\end{pmatrix}.$\end{proof}

\subsection{Polar decompositions} A specially remarkable instance occurs when one sets the polar decomposition $T=UP=U(T^*T)^{1/2}$ of an operator $T\in \mathfrak L(\ell_2) $.

\begin{prop}
Let $T\in \mathfrak L(\ell_2)$ be an operator and $T=UP=U(T^*T)^{1/2}$ its polar decomposition. If $U$ is an operator on the scale then
$\begin{pmatrix}
U & T\\
O & U
\end{pmatrix}$ is a symplectic transformation on $Z_2$.
\end{prop}
\begin{proof}

By Proposition \ref{opT} we just have to recall that $U^*T = P$ is selfadjoint.\end{proof}

Thus, for every selfadjoint operator $T$ the operator $\begin{pmatrix}
U & UT\\
0 & U
\end{pmatrix}$ is a symplectic transformation.

\subsection{Diagonal operators}\label{Subsubsection_03}
Let $\sigma \in \ell_\infty$. The diagonal operator $\sigma\big((x_n)_n\big)=(\sigma_nx_n)_n$ is an operator on the scale and it therefore induces the operator $\tau_\sigma = \begin{pmatrix}
\sigma & 0\\
0 & \sigma
\end{pmatrix}$ on $Z_2$. The operator $\tau_\sigma $ is a symplectic transformation if and only if $\overline{\sigma_n}\sigma_n=1$ for all $n\geq1$ since $\sigma^* = \overline \sigma$ and

$$\tau_\sigma^+\tau_\sigma = \begin{pmatrix}
\sigma & 0\\
0& \sigma
\end{pmatrix}^+\begin{pmatrix}
\sigma & 0\\
0& \sigma
\end{pmatrix}=\begin{pmatrix}
\sigma^* & 0\\
0 & \sigma^*
\end{pmatrix}\begin{pmatrix}
\sigma & 0\\
0& \sigma
\end{pmatrix}=\begin{pmatrix}
1 & 0\\
0 & 1
\end{pmatrix}.$$

Thus, $\tau_\sigma$ is a selfadjoint symplectic transformation if and only if $\sigma \in \{-1,1\}^\N$. 

\subsection{Shift operators}\label{Subsubsection_02} The right-shift operator $r((x_n)_n)=\big((x_{n-1})_n\big)$ is an isometric operator on the scale and therefore $R = \begin{pmatrix}
r & 0\\
0 & r
\end{pmatrix}\in \mathfrak L(Z_2)$ is an isometry on $Z_2$ with $2$-codimensional range.
The adjoint $\ell=r^*$ is the left-shift operator $\ell((x_n)_n)=\big((x_{n+1})_n\big)$ is also a operator on the scale and therefore $L=\begin{pmatrix}
\ell & 0\\
0 & \ell
\end{pmatrix} \in\mathfrak L(Z_2)$. It follows from Proposition \ref{opT} that $R$ is a symplectic transformation (see also below), while  $L$ is not symplectic because it is not injective. The comments at the end of the previous section imply that $R^+=L$ and thus $LR=R^+R=I$.

\subsection{Block operators} \label{Block operators}Let $\mathfrak u$ be a sequence $(u_n)_n$ of disjointly supported normalized blocks in $\ell_2$, that we can understand as the operator $\mathfrak u: \ell_2\to \ell_2$ given by $\mathfrak u(x)= \sum x_n u_n$. In general $\mathfrak u$ is not an operator on the scale and $\begin{pmatrix}
\mathfrak u  & 0\\
0 & \mathfrak u
\end{pmatrix}$ is not an operator in $Z_2$. The \emph{block operator} $T_\mathfrak u :Z_2\rightarrow Z_2$ is defined as $T_\mathfrak u(e_n, 0)=(u_n,0)\quad\text{and}\quad T_\mathfrak u(0,e_n)=(\KP u_n, u_n)$, namely
\begin{equation}\label{eq_06} T_\mathfrak u = \begin{pmatrix}
\mathfrak u  & \KP \mathfrak u\\
0 & \mathfrak u
\end{pmatrix}.
\end{equation}

The operators $T_\mathfrak u$ are symplectic transformations, and the proof for this can be followed in detail in \cite[Section 10.9: The Properties of $Z_2$ explained by itself]{CabeCast2021}. We will prove the general case in Section \ref{blocksymp}. The idea is that equation (\ref{eq_05}) is equivalent to $T_{\mathfrak u}^+T_{\mathfrak u}=I$ or else to  $D=T_{\mathfrak u}^*DT_{\mathfrak u}$, where $D:Z_2\rightarrow Z_2^*$ is the duality isomorphism given in Proposition \ref{P_01} and this follows, after a few cumbersome computations, from the equality
$$\left \langle u_j,\sum_ix_i \KP(u_i) \right \rangle = \left \langle \KP(u_j),\sum_ix_iu_i \right\rangle.$$
The operator $T_{\mathfrak u}D^{-1}T_{\mathfrak u}^*D$ therefore defines a projection onto $T_{\mathfrak u}[Z_2]$ and this shows that $T_{\mathfrak u}$ is an into isometry (respect to the usual quasi-norm) with complemented range. All these results are from \cite{kaltsym}.

Observe that $R$ can be regarded as a block operator with the choice of $(u_n)_{n\in\mathbb{N}}=(e_{n+1})_{n\in\mathbb{N}}$ since $\KP e_n=0$.

\subsection{Transvections} Throughout this section $(X,\omega)$ will denote a symplectic Banach space, in particular $Z_2$. The \emph{symplectic group} of  $(X,\omega)$  is the subgroup $\Sp(X,\omega)$ of $\GL(X)$ of all symplectic automorphisms:
$$\Sp(X, \omega)=\{T \in \GL(X) \colon \omega(Tx,Ty) = \omega (x,y) \text{ for all }x,y\in X\}.$$ We will denote by  $\Sp (Z_2)$ the symplectic group of $Z_2$ endowed with the  symplectic form of Prop. \ref{P_01}. Observe that $\Sp (Z_2)$ is not a bounded subgroup. Indeed, if $D_a:\ell_2\rightarrow\ell_2$ is a diagonal operator on $\ell_2$ given by some real $a\in\ell_\infty$ then $D_a$ is selfadjoint and $\begin{pmatrix}
	I & D_a\\
	0 & I
	\end{pmatrix}$ is an invertible symplectic transformation with norm  $\|a\|_\infty+1$.

\adef Let $\lambda\in\mathbb{K}$ and $u\in X$. The \emph{transvection} associated to $\lambda$ and $u$ is the linear map $\mathcal{T}_{u,\lambda}$ given by $\mathcal{T}_{u,\lambda}(x)=x+\lambda \omega( x,u ) u,\quad\text{for each }x\in X$.\zdef

Recall that for a subset $U\subset X$, we denote by $U^{\omega}=\{x\in X\colon \omega( x,u ) =0,\text{ for all }u\in U\}$ the symplectic anihilator of $U$. For any $u\in X$, let us denote by $u^{\omega}$ the  anihilator of $\{u\}$. By duality of Lemma \ref{subsym}, the annihilator $u^{\omega}$ of a line defined by any $u\in X$ is an hyperplane of $X$. It follows that a transvection $\mathcal{T}_{u,\lambda}$ is the identity on the hyperplane $u^{\omega}$ and the identity on the corresponding quotient $X/u^{\omega}$ (as $\mathcal{T}_{u,\lambda}(x)-x\in u^{\omega}$ for each $x\in X$).

\begin{lemma}\label{L_transveccion} Transvections are symplectic transformations.
\end{lemma}
\begin{proof}
In the first place, transvections are linear due as the symplectic form $\omega( \cdot,\cdot )$ is bilinear (sesquilinear in the complex case). Given $x\in X$, boundedness follows by
\begin{align*}
\|\mathcal{T}_{u,\lambda}(x)\|_{X}&\leq \|x\|_{X}+|\lambda|\,|\omega( x,u )|\|u\|_{X}\\
&\leq\big(1+|\lambda|\,\|\omega\|\,\|u\|_{X}^2\big)\|x\|_{X}.
\end{align*}
Now taking $x,y\in X$ we deduce that
\begin{align*}
\omega( \mathcal{T}_{u,\lambda}(x),\mathcal{T}_{u,\lambda}(y) ) &= \omega( x+\lambda \omega( x,u ) u, y+\lambda \omega( y,u ) u )\\
&=\omega( x,y ) + \lambda \omega( y,u )\omega( x,u )-\lambda \omega( x,u )\omega(y,u)\\
&=\omega( x,y ).\qedhere
\end{align*}
\end{proof}

Moreover, observe that  $\mathcal{T}_{u,\lambda}\mathcal{T}_{u,\mu}=\mathcal{T}_{u,\lambda+\mu}$ and $\mathcal{T}_{au,\lambda}=\mathcal{T}_{u,a^2\lambda}$ as it immediately follows from the definition of transvection, which shows that, given $u\in X$, $\mathscr{T}_u=\{\mathcal{T}_{u,\lambda}\colon \lambda \in \mathbb K\}$ is a subgroup of $\Sp (X)$ and the map $\lambda\in(\mathbb{K},+)\mapsto\mathcal{T}_{u,\lambda}$ defines an isomorphism of groups.

\section{Rochberg spaces are symplectic}

Consider the complex interpolation method applied to the scale $(\ell_\infty, \ell_1)$ (see the classical \cite{belo}; or else \cite{CabeCast2021}). It is well known that it provides the space $(\ell_\infty, \ell_1)_{\theta} = \ell_{\theta^{-1}}$ for $0<\theta<1$. In particular, $(\ell_\infty, \ell_1)_{1/2} = \ell_{2}$. The Rochberg spaces \cite{rochberg} obtained at $1/2$ are defined as
$$\mathfrak R^{(n)}=\{(x_{n-1},\ldots,x_1,x_0)\in\ell_\infty^n\colon x_i=f^{(i)}(1/2)/i!,\text{ for some } f\in\mathcal{C},\,0\leq i\leq n-1\}$$
here $\mathcal C$ represents the associated Calder\'on space) and can be considered as generalizations of $\ell_2$. Indeed,  $\mathfrak{R}^{(1)}=\ell_2$ and it was Kalton who noticed that $\mathfrak{R}^{(2)}=Z_2$ (see  \cite{Pacific,tams} for additional information). To show that Rochberg spaces are symplectic we need first to know that they are isomorphic to their duals in the following form taken from \cite{CabCasCorr}:

\begin{prop}\label{dualRochberg}
Consider for each $n\geq1$ the continuous bilinear map $\omega_n:\mathfrak R^{(n)}\times\mathfrak R^{(n)}\rightarrow\mathbb{R}$ given by
$$\omega_n\big((x_{n-1},\ldots,x_0),(y_{n-1},\ldots,y_0)\big)=\sum_{i+j=n-1}(-1)^i\langle x_i,y_j\rangle.$$
The induced operator $D_n:\mathfrak R^{(n)}\rightarrow\mathfrak{R}^{{(n)}^*}$ given by $D_n(x)(y)=\omega_n(x,y)$ is an isomorphism onto.
\end{prop}

That this duality makes Rochberg spaces symplectic for even $n$, as it occurs with $Z_2$ is, somehow, expected.
Surprisingly enough, odd Rochberg spaces are also symplectic, but not in the same way as even Rochberg spaces. In fact, in the Hilbert space case we see that there is a correspondence between complex and symplectic structures: if $\omega$ is a symplectic form, then there is a complex structure $J$ such that $\omega(x,y)=\langle x,J(y) \rangle$. Namely, a symplectic structure is obtained ``twisting'' the natural duality with a complex structure. This approach generalizes to higher odd Rochber spaces, i.e., a complex structure on $\Rn$ may be used to induce a perturbation on $\omega_n$ and define a symplectic structure.

\begin{theorem}\label{RochbergSymplectic} All Rochberg derived spaces are symplectic.
\end{theorem}
\begin{proof} Observe that $\omega_n$ is alternated if and only if $n$ is even; so, the result holds for even $n$.\medskip

For $n$ odd, consider a complex structure $\sigma$ on $\ell_2$ that is an operator on the scale; say, $\sigma(x)=(-x_2,x_1,-x_4,x_3,\ldots)$, so that the induced diagonal operator $\tau_\sigma$ is bounded on $Z_2$. The generalized form of the commutator theorem, see \cite{racsam,CabCasCorr}, shows that the $n\times n$ matrix diagonal operator still acts boundedly on the corresponding $\Rn$. We will continue calling $\tau_\sigma$ this diagonal operator.  We define the bilinear map

$$\overline{\omega_{n}}\big((x_{n-1},\ldots,x_0),(y_{n-1},\ldots,y_0)\big)= \omega_{n}\big((x_{n-1},\ldots,x_0),\tau_\sigma(y_{n-1},\ldots,y_0)\big)=\sum_{i+j=n-1}(-1)^{i}\langle x_i,\sigma y_j \rangle.$$

This map is now alternated due to the fact that $\sigma^*=-\sigma$. Indeed,
\begin{align*}
\overline{\omega_{n}}\big((x_{n-1},\ldots,x_0),(y_{n-1},\ldots,y_0)\big)&=\sum_{i+j=n-1}(-1)^{i}\langle x_i,\sigma y_j \rangle\\
&=\sum_{i+j=n-1}(-1)^{i}\langle \sigma^*x_i, y_j \rangle\\
&=\sum_{i+j=n-1}(-1)^{i}\,(-1)\langle \sigma x_i, y_j \rangle\\
&=(-1)\,\sum_{i+j=n-1}(-1)^{i}\langle y_j, \sigma x_i \rangle\\
&=(-1)\,\sum_{j+i=n-1}(-1)^{i}\,(-1)^{i+j}\langle y_j, \sigma x_i \rangle\\
&=(-1)\,\sum_{i+j=n-1}(-1)^{j}\langle y_j, \sigma x_i \rangle\\
&=-\overline{\omega_{n}}\big((y_{n-1},\ldots,y_0),(x_{n-1},\ldots,x_0)\big).
\end{align*}

Boundedness follows from the boundness of $\omega_n$ and $\tau_\sigma$:
$$\big|\overline{\omega_{n}}(x,y)\big|=\big|\omega_{n}(x,\tau_\sigma y)\big|\leq K\,\|x\|\,\|\tau_\sigma y\|\leq C \|x\|\,\|y\|.$$

To obtain that $(\Rn,\overline{\omega_n})$ is symplectic it suffices to show that the induced linear map $L_{\overline{\omega_{n}}}:\Rn\rightarrow{\Rn}^*$ is an isomorphism onto. Assume that there exists $x\in\Rn$ such that $L_{\overline{\omega_{n}}}(x)(y)=0$ for all $y\in\Rn$. Thus $L_{\omega_n}(x)(\tau_\sigma y)=0$ for all $y\in\Rn$. Taking into account that $\tau_\sigma$ is invertible in $\Rn$, it follows that $L_{\omega_n}(x)(y)=0$ for all $y\in\Rn$, so that $x=0$. Moreover, as $\tau_\sigma$ is an isomorphism, its clear that $\overline{\omega_n}$ has closed range because $\omega_n$ has it.\end{proof}

\section{Block operators on Rochberg spaces are symplectic}\label{blocksymp}

A sequence $\mathfrak u= (u_n)_{n\in\mathbb{N}}$ of normalized blocks in $\ell_2$ induces a multiplication operator $u:\ell_2\rightarrow\ell_2$ given by $u(e_n)=u_n$ and, as we showed in section \ref{Block operators}, a block operator
$\begin{pmatrix}
\mathfrak u & 2\mathfrak u\log \mathfrak u \\
0 & \mathfrak u &\\
\end{pmatrix}$ in $Z_2$. Their higher order generalizations of block operators were obtained in \cite[Proposition 7.1]{CasFer21} as the operators
$$T_{\mathfrak u, n} =\begin{pmatrix}
\mathfrak u & 2\mathfrak u\log \mathfrak u & 2\mathfrak u\log^2\mathfrak u & \cdots & \frac{2^{n-1}}{(n-1)!}\mathfrak u\log^{n-1}\mathfrak u\\
0 & \mathfrak u & 2\mathfrak u\log \mathfrak u & 2\mathfrak u\log^2\mathfrak u  & \cdots\\
0 & 0 & \mathfrak u & 2\mathfrak u\log \mathfrak u & 2\mathfrak u\log^2\mathfrak u\\
0 & 0 & 0 & \mathfrak u & 2\mathfrak u\log \mathfrak u\\
0 & 0 & 0 & 0 & \mathfrak u &\\
\end{pmatrix}  $$
that act boundedly $T_{\mathfrak u, n}:\mathfrak R^{(n)}\rightarrow\mathfrak R^{(n)}$. See also \cite[Section 5.1]{Pacific}. We will shorten the name to $T_{\mathfrak u}$ when no confusion is possible about which is $n$. Given $T\in \mathfrak L(\mathfrak R^{(n)})$, recall that $T^+$ always denotes the symplectic adjoint of an operator $T$, namely $\omega_n(T^+x, y) = \omega_n(x, Ty)$. If $T\in \mathfrak L(\mathcal R^{(2n-1)})$ then we will denote
$T^\sharp$ the symplectic adjoint with respect to $\overline \omega_n$, namely $\overline \omega_n(T^\sharp x, y) = \overline \omega_n(x, Ty)$. We have:

 \begin{lemma}\label{P_06}
	Let $T=\begin{pmatrix}
	A_0 & B_0 & C_0 & \cdots & D\\
	0 & A_1 & B_1 & C_1  & \cdots\\
	0 & 0 & A_2 & B_2 & C_{n-3}\\
	0 & 0 & 0 & A_3 & B_{n-2}\\
	0 & 0 & 0 & 0 & A_{n-1} &\\
	\end{pmatrix}$ be an upper triangular operator on $\mathfrak R^{(n)}$.
\begin{itemize}
\item $T^+=\begin{pmatrix}
	A_{n-1}^* & -B_{n-2}^* & C_{n-3}^* & \cdots & (-1)^{n-1}D^*\\
	0 & A_{n-2}^* & -B_{n-3}^* & C_{n-4}^*  & \cdots\\
	0 & 0 & A_{n-3}^* & -B_{n-4}^* & C_0^*\\
	0 & 0 & 0 & A_{n-4}^* & -B_0^*\\
	0 & 0 & 0 & 0 & A_0^* &\\
	\end{pmatrix}.$
\item If $n$ is odd, $T^\sharp = - T^+ \tau_\sigma$.\end{itemize} \end{lemma}
\begin{proof} The first part can be obtained by plain induction. The second part is simple: since $ \overline \omega_n(x, Ty) = \omega_n(x, \tau_\sigma T y)$ then
$$T^\sharp = (\tau_\sigma T)^+ = T^+ \tau_\sigma^+ = T^+ (\tau_\sigma)^* = - T^+ \tau_\sigma.\qedhere$$
\end{proof}

We prove now that block operators are symplectic.

\begin{prop}\label{P_08} Let $D_n:\mathfrak R^{(n)}\rightarrow\mathfrak{R}^{{(n)}^*}$ be the duality isomorphism $D_n(x)(y)=\omega_n(x,y)$ from Proposition \ref{dualRochberg}. One has $T_{\mathfrak u}^*D_n T_{\mathfrak u} = D_n$ or, equivalently, $\omega_n(T_{\mathfrak u} x, T_{\mathfrak u} y)=\omega_n(x, y).$
\end{prop}
\begin{proof} Reduced to its bare bones the argumentation says that if we call $x_{i,k}$ to the vector having $e_i$ at the $k^{th}$ position and $k+l=n+1$ then
\begin{align}\label{eq_08}
\omega_n (T_{\mathfrak u} (x_{i,k})&,T_{\mathfrak u}(x_{j,l})\\
&=\omega_n\Big(\big(\frac{2^{k-1}}{(k-1)!}u_i\log^{k-1}|u_i|,\ldots,u_i,\ldots,0\big),\big(\frac{2^{l-1}}{(l-1)!}u_j\log^{l-1}|u_j|,\ldots,u_j,\ldots,0\big)\Big)\nonumber\\
&=(-1)^k\langle u_i,u_j \rangle=\omega_n(x_{i,k},x_{j,l})\nonumber.
\end{align}

If $k+l<n+1$ then (\ref{eq_08}) cancels out as we are multiplying by zeroes. If $k+l>n+1$ then (\ref{eq_08}) becomes, after setting  $m=k+l-n$,
\begin{align*}
(-1)^k \Big\langle u_i,\frac{2^m}{m!}u_j\log^{m}|u_j|\Big\rangle&+(-1)^{k-1} \Big\langle \frac{2^1}{1!}u_i\log|u_i|,\frac{2^{m-1}}{(m-1)!}u_j\log^{m-1}|u_j| \Big\rangle\\
&+\cdots+(-1)^{k-m}\Big\langle \frac{2^m}{m!}u_i\log^{m}|u_i|,u_j \Big\rangle.
\end{align*}

If $i\neq j$ then this last summand is null because $\langle u_i,u_j \rangle=0$. If $i=j$ then (\ref{eq_08}) becomes
\begin{align*}
\log^m|u_i|\Big[\sum_{p=0}^m(-1)^{k+p}\frac{2^m}{p!(m-p)!}\Big]&=\log^m|u_i|\Big[\sum_{p=0}^m(-1)^{k+p}\frac{2^m}{m!}{m\choose p}\Big]\\
&=\log^m|u_i|\,\frac{2^m}{m!}(-1)^k\Bigg[\sum_{p=0}^m(-1)^p{m\choose p} \Bigg].
\end{align*}
Now, the Binomial Theorem $0=(1-1)^m=\sum_{k=0}^m{m\choose k}1^{m-k}(-1)^k$ cancels out all terms.\end{proof}

We conclude this section with several technical lemmata of independent interest about generalized block operators.

\begin{lemma}\label{complementedrange} The range $T_{\mathfrak u, n}[\mathfrak R^{(n)}]$ is isomorphic to $\mathfrak R^{(n)}$ and complemented in $\mathfrak R^{(n)}$.\end{lemma}
\begin{proof} The proof follows the arguments of \cite{KalSwan,symmetries21}. There is a commutative diagram

$$\begin{CD}
0@>>>  \mathfrak R_{n-1}@>>> \mathfrak R_n @>>> \ell_2@>>> 0\\
 &&@VT_{\mathfrak u, n-1}VV @VVT_{\mathfrak u, n}V @VVT_{\mathfrak u, 1}V \\
0@>>>  \mathfrak R_{n-1}@>>> \mathfrak R_n @>>> \ell_2@>>> 0
\end{CD}$$
By a general $3$-space property (see \cite{CabeCast2021}), since $T_{\mathfrak u, 1} = \tau_{\mathfrak u}$ is an isometry, $T_{\mathfrak u, n}$ must be into isometries. Thus, $T_{\mathfrak u, n} [\mathfrak R^{(n)}]$ is an isometric copy of $\mathfrak R^{(n)}$. This isometric copy is complemented because of the identity $(T_{\mathfrak u, n})^* D_n T_{\mathfrak u, n} =D_n$ from Proposition \ref{P_08}. Thus, $T_{\mathfrak u, n} D_{n}^{-1}(T_{\mathfrak u, n})^*D_{n}$ is a projection of $\mathfrak R^{(n)}$ onto the range of $T_{\mathfrak u, n}$. \end{proof}

We now extend the classical result about the behaviour of operators on $Z_2$ due to Kalton \cite[Lemma 6]{kaltsym} to higher order Rochberg spaces.

\begin{lemma} If $T:\mathfrak R^{(n)}\rightarrow\mathfrak R^{(n)}$ is not strictly singular then there exists $\alpha\neq0$ and block operators $T_{\mathfrak u}$ and $T_{\mathfrak v}$ such that $TT_{\mathfrak u}-\alpha T_{\mathfrak v}$ is strictly singular.\end{lemma}
	\begin{proof} Let us recall from \cite{Pacific} that the canonical exact sequence \begin{eqnarray} \label{ss} \xymatrix{0\ar[r]& \ell_2 \ar[r] &\mathfrak R^{(n)}\ar[r] &\mathfrak{R}^{(n-1)}\ar[r] & 0}\end{eqnarray} has strictly singular quotient map. Therefore, there exists $\alpha\neq0$ and normalized block basic sequences $\mathfrak u = (u_n)_n$ and $\mathfrak v= (v_n)_n$ in $\ell_2$ such that
		$$T(u_n,0,0,\ldots,0)=\alpha(v_n,0,0,\ldots,0)+w_n,\quad\text{for each }n\in\mathbb{N},$$
		where $\sum_n\|w_n\|<\infty$. Now just take the block operators $T_{\mathfrak u},T_{\mathfrak v}$ induced by those sequences and define $K(e_n,0,0,\ldots,0)=w_n$ for all $n\in \mathbb{N}$ and $K(0,0,\ldots,\overset{(i)}{e_n},\ldots,0)=0$ for all $n\in\mathbb{N}$ and $1\leq i\leq n-1$. This $K$ is compact and $TT_{\mathfrak u}-\alpha T_{\mathfrak v}-K=0$ on the canonical copy of $\ell_2$ in $\mathfrak R^{(n)}$ given by (\ref{ss}) and thus  $TT_{\mathfrak u}-\alpha T_{\mathfrak v}$ must be strictly singular.\end{proof}

An important consequence of this perturbation result is:
\begin{lemma} Let $T\in \mathfrak L \left( \mathfrak R^{(n)}\right) $. \begin{itemize}
\item If $T^+T$ is strictly singular then $T$ is strictly singular.
\item If $n$ is odd and $T^\sharp T$ is strictly singular then $T$ is strictly singular.
\end{itemize}
\end{lemma}
\begin{proof} In general, if $T$ is not strictly singular there exists $\alpha\neq0$ and block operators $T_{\mathfrak u},T_{\mathfrak v}$ such that $ TT_{\mathfrak u}=\alpha T_{\mathfrak v}-S$ with $S$ strictly singular. Therefore $$T_{\mathfrak u}^+T^+TT_{\mathfrak u}=(TT_{\mathfrak u})^+TT_{\mathfrak u}=(\alpha T_{\mathfrak v}^+-S^+)(\alpha T_{\mathfrak v}-S)=\alpha'T_{\mathfrak v}^+T_{\mathfrak v}+S'=\alpha'I+S'$$ and $I$ must be strictly singular. If $n$ is odd, the last calculation becomes
$$T_{\mathfrak u}^\sharp T^\sharp TT_{\mathfrak u}=(TT_{\mathfrak u})^\sharp T T_{\mathfrak u}=(\alpha T_{\mathfrak v}^\sharp-S^\sharp)(\alpha T_{\mathfrak v}-S)=\alpha'T_{\mathfrak v}^\sharp T_{\mathfrak v} + S'=
-\alpha'T_{\mathfrak v}^+\tau_\sigma T_{\mathfrak v} + S'$$
by Proposition \ref{P_06}. This means that $T_{\mathfrak v}^+\tau_\sigma T_{\mathfrak v}$ must be strictly singular, but since $\tau_\sigma T_{\mathfrak v}$ is invertible,
$T_{\mathfrak v}^+$ must be strictly singular, as well as $T_{\mathfrak v}$, which therefore must be finite dimensional, hence $TT_{\mathfrak u}$ must be strictly singular as well as $T$.\end{proof}

\section{Rochberg spaces do not contain Lagrangian subspaces}

We now extend the Kalton-Swanson theorem \cite{KalSwan} showing that the symplectic structures of Rochberg spaces we have defined are not trivial.

\begin{theorem}\label{nol2}	$\mathfrak{R}^{(n)}$ has no Lagrangian subspace.\end{theorem}
\begin{proof} Let $T$ be a projection onto an infinite dimensional isotropic subspace. If $n$ is even, $T^+T=0$, so $T$ must be strictly singular and thus every complemented isotropic subspace must be finite dimensional. If $n$ is odd, $T^\sharp T=0$ and then $T$ must be strictly singular. Therefore every complemented isotropic subspace must be finite dimensional.\end{proof}

Some of the authors of this paper conjecture that Rochberg spaces obtained from a reflexive Banach space $X$ such that $X \cap \overline X^*$ is dense in both $X$ and $\overline X^*$ and $(X, \overline X^*)_{1/2}$ is isometric to a Hilbert space are symplectic. Conditions to obtain  $(X, \overline X^*)_{1/2}$ are in  \cite{Wat} (see also \cite{tams}). By \cite{CabCasCorr}, see also \cite[Proposition 2.11]{rochweiss}, the diagram
\begin{equation}\label{duality}\begin{CD}
0 @>>> \ell_2^* @>>>  (\mathfrak{R}^{(2)})^* @>>> \ell_2^*@>>>0  \\
@. @|  @VV{T}V @| \\
0@>>> \ell_2^*@>>> \mathfrak{R}^{*(2)}  @>>> \ell_2^*@>>>0
\end{CD}\end{equation}
with $T(x^*,y^*)(x,y)= \pin{x^*}{y}+ \pin{y^*}{x}$ for all $(x^*,y^*)\in d(X_\theta^*)$ and all $(x,y)\in dX_\theta$  is commutative, which shows that the corresponding $\mathfrak R^{(2)}$ is symplectic. The assertion generalizes to ``even Rochberg spaces $\mathfrak{R}^{(2n)}$ are symplectic". Moreover, if both $X$ and $X^*$ have a common complex structure then, reasoning as in Theorem \ref{RochbergSymplectic}, we deduce that also odd Rochberg spaces $\mathfrak{R}^{(2n-1)}$ are symplectic. However, it seems difficult to prove that such symplectic structures are non-trivial. If one wants to mimicry the proof for the couple $(\ell_1,\ell_\infty)$ to the general case one should prove that the corresponding block operators
$$\left(
    \begin{array}{ccccc}
      u & \Omega^{(1)}(u) & \Omega^{(2)}(u) & \dots & \Omega^{(n-1)}(u) \\
      0 & u & \Omega^{(1)}(u)  & \Omega^{(2)}(u) & \dots \\
      0 & 0 & u & \Omega^{(1)}(u)  & \Omega^{(2)}(u) \\
      \dots & 0 & 0 & u & \Omega^{(1)}(u) \\
      0 & \dots  & 0 & 0 & u \\
    \end{array}
  \right)$$
where $\Omega^{(n)}$ are the corresponding $n$-differentials (see \cite{rochberg})  are symplectic operators. See also \cite{CasFer21} to determine how they can be calculated. In our Kalton-Peck case, $\Omega^{(n)} x = \frac{2^n}{n!}x\log^n x$ for normalized $x$. In general, it is perfectly possible (think about the case of weighted Hilbert spaces \cite{symmetries21}) that all iterated Rochberg spaces are actually Hilbert spaces and thus their symplectic structures are trivial.

\section{Almost symplectic structures} \label{almost}

This section is motivated by the study \cite{ccfm} on complex structures on $Z_2$. Let us recall that a linear bounded operator $J:X\to X$ defined on a real Banach space is a \emph{complex structure} if $J^2=-I$; and that a complex structure $J$ on $X$ yields a $\C$-linear structure on $X$ by declaring $(\alpha + i\beta)x= \alpha x + \beta J(x).$ The resulting complex space will be denoted $X^J$ and it is a Banach space equipped with the norm $|||x||| = \sup_{\theta\in[0,2\pi]} \| \cos \theta x + \sin \theta J(x)\|$.

The following result of Weinstein \cite{wein} shows that symplectic structures on real Hilbert spaces are obtained from complex structures.
\begin{lemma} Let $\mathcal H$ be a real Hilbert space. For every symplectic form $\omega$ in $\mathcal H$ there exist a complex structure $J$ on $\mathcal H$ and an equivalent inner product $\pin{\,}{\,}_R$ on $\mathcal H$  such that $\omega(x,y)=\pin{x}{Jy}_R$ for every $x,y\in \mathcal H$. \end{lemma}

Consequently, all symplectic structures on a real Hilbert space are equivalent to the standard one. Let us prove that a complex structure on a real Hilbert space induces a symplectic structure.

\begin{lemma}  Let $\mathcal H$ be a real Hilbert space and $J$ be a complex structure on $\mathcal H$. Then there exist a symplectic form $\omega$ on $\mathcal H$ and an equivalent inner product $\pin{\,}{\,}_R$ on $\mathcal H$ such that $\omega(x,y)=\pin{x}{Jy}_R$ for every $x,y\in \mathcal H$.
\end{lemma}
\begin{proof} Let  $J$ be a complex structure on $\mathcal H$. Let us take $R=I+J^{*}J$. Then $\pin{x}{y}_R=\pin{x}{Ry}$ defines  an  equivalent inner product on $\mathcal H$ for which $J$ is an isometry and therefore an unitary operator. The Hilbert-adjoint of $J$ with respect to this inner product is $J^{-1}=-J$. Then $\omega(x,y):=\pin{x}{Jy}_R$ is symplectic on $\mathcal H$. \end{proof}


It is straightforward that a complex structure on a hyperplane of any Banach space cannot be extended to a complex structure on the whole space. We have that the same situation occurs for symplectic structures (Corollary  \ref{symio}).  In \cite{ccfm} we showed that no complex structure on $\ell_2$ can be extended to a complex structure on a hyperplane of $Z_2$ containing it. We now observe that an analogous result holds for symplectic structures (Corollary \ref{symhyp}).

\begin{definition}
Let $X$ and $Y$ be Banach spaces and $j: Y \to X$ be an isomorphism into. A bilinear map $\Omega$ on $X$ extends a bilinear map $\omega$ on $Y$ \emph{through} $j$ when $\Omega(jx,jy)=\omega(x,y)$ for every $x,y\in Y$.
Equivalently  $\Omega$ extend $\omega$ through $j$ if the diagram is commutative
\begin{equation}\label{extension}\begin{CD}
 Y @>j>> X  \\
@V{L_\omega}VV @VV{L_\Omega}V \\
Y^* @<j^*<< X^*
\end{CD}\end{equation} \end{definition}

Essentially following \cite[Prop. 3.1]{cgj} we obtain

\begin{prop} Let $X$ and $Y$ be Banach spaces and $j: Y \to X$ be an isomorphism into. If a symplectic structure $\omega$ on $Y$ can be extended to a bilinear form $\Omega$ on $X$, then $j(Y)$ is complemented on $X$.
\end{prop}
\begin{proof} Indeed, $jL_\omega^{-1}j^*L_\Omega$ would be a projection onto $j[Y]$ for symplectic extensions.
\end{proof}
Since $Z_2$ does not contain complemented copies of $\ell_2$ \cite[Corollary 6.7]{Kalton_Peck} we have
\begin{cor} \label{symhyp}
No symplectic structure on $\ell_2$ can be extended to a bilinear form on a hyperplane $H$ of   $Z_2$  through any embedding $j: \ell_2 \to H$.
\end{cor}

We  observe now that a symplectic structure on an hyperplane induces an almost symplectic structure on the space in the following sense:
\begin{definition} Let $X$ be a  real (complex, resp.) Banach space and  let $\alpha:X\to X^*$ ($\alpha: X \to \overline X^*$, resp.) be an isomorphism.  We say that $\alpha$ is almost symplectic  if $\alpha +\alpha^*$ ($\alpha + \overline \alpha^*$, resp.) is strictly singular. We will say that $X$ admits an almost symplectic structure
if there exists an almost symplectic isomorphism $X\to X^*$ ($X \to \overline X^*$, resp.).\end{definition}

Recall that for a given Banach space $X$ and $F$ a subspace of $X$ the annihilator of $F$ is the closed subspace of $X^*$ defined by $F^\perp=\{ f\in X^* \, : \, f(x)=0 \, \, \text{for all} \, \, x\in F\}$.

Now, if $\beta$ is a symplectic structure on a hyperplane $H$  of $X$  and we identify $H^*=X^*/H^{\perp}$, then we can consider the isomorphism extension $\alpha : H\oplus[e] \to H^{\perp} \oplus [e]^{\perp}$ given by $\alpha(h+\lambda e)= \beta(h)+\lambda e^*$. Quite clearly $\alpha + \alpha ^*$ is a rank one operator. We prove now the converse:

\begin{prop}\label{almostsym} Let $X$ be a  real  or complex Banach space admitting an almost symplectic structure. Then either $X$ or its hyperplanes admit a symplectic structure. In particular, $X$ is reflexive.
\end{prop}

\begin{proof}  In the real case, consider  $\alpha:X\to X^*$ an almost symplectic isomorphism and let $s:X\to X^*$ be a strictly singular operator such that $\alpha + \alpha^*= s$. Then, denoting by $\beta= \alpha -s/2$,  we have $ \beta^*=-\beta$ and that $\beta$ is a Fredholm operator with index 0. By Fredholm theory there exist closed subspaces $X_0\subseteq X$ and $Y_0\subseteq X^*$ such that $E= X_0 \oplus \ker \beta$ and $X^*= Y_0 \oplus F$, where  $\ker \beta$ and $F$ are finite dimensional subspaces with the same dimension and such that the restriction $\gamma:= \beta_{|X_0} \, : X_0 \to Y_0$ is an isomorphism onto.

Observe that that $Y_0=(\ker \beta) ^{\perp}$. Indeed, let $\phi\in Y_0$  and $x\in \ker \beta$. Let $x_0\in X_0$ be such that $\phi=\beta(x_0)$, then $\phi(x)=\beta(x_0)(x)=-\beta(x)(x_0)=0$. Hence $Y_0\subseteq (\ker \beta) ^{\perp}$. The equality holds since $\di (X^*/(\ker \beta)^{\perp})= \di (X^*/Y_0)$. Then we can assume that $F= X_0^\perp$. Let us take $\{x_1, \ldots, x_n\}$ and $\{\phi_1, \ldots , \phi_n\}$ basis of $\ker \beta$ and $X_0^\perp$, respectively, such that $\phi_i(x_j)= \delta_{ij}$ for all $1\leq i,j\leq n$. Suppose now that $n$ is even and consider the map  $\psi \, : \, \ker \beta \to X_0^\perp$ defined by  $\psi(x_{2k-1})= \phi_{2k}$ and  $\psi(x_{2k})=-\phi_{2k-1}$ for $k=1, \ldots , n/2$.  It follows that the map  $\Gamma \, : \, X_0 \oplus \ker \beta \to (\ker \beta) ^{\perp} \oplus X_0^\perp$ defined by the matrix
$\begin{pmatrix} \displaystyle
 \gamma & 0\\
  0 & \psi
\end{pmatrix}$  is symplectic  on $X$.
When $n$ is odd, the previous construction gives us a symplectic structure for an hyperplane of $X$. The proof for the complex case is completely analogous.\end{proof}

This proof implies that if $X$ is a Banach space (over $\R$ or $\C$) and $\alpha \, : \,X\to X^*$ ($\alpha \, : \,X\to \overline{X}^*$, resp.) is an isomorphism such that $\alpha - \alpha^*$  ($\alpha - \overline{\alpha}^*$, resp.) is strictly singular then there exists an isomorphism $\beta \, :  \,X \to X^*$ ($\beta \, :  \,X \to \overline{X}^*$, resp.) such that $\beta^* =\beta $ ($\overline {\beta}^* =\beta $, resp.). In both  the real and complex case, we call such isomorphism $\beta$  a \emph{Hermitian structure} on $X$.


\begin{prop}  Let  $X$  be a reflexive Banach space over $\mathbb K $ ($\R$ or $\C$)  isomorphic to its dual (dual conjugate when $\mathbb K= \mathbb C$) such that every operator $T \, : \,X\to X$ is of the form $\lambda I +S$  for some $\lambda \in \mathbb K$ and $S$ strictly singular. Then we have either an induced almost symplectic structure or a Hermitian structure on $X$. \end{prop}

\begin{proof} For the real case,  let $\alpha$ be an isomorphism onto the dual. Using the canonical identification of $(\alpha) ^{-1}\alpha^* \, : \, X^{**} \to X$ as an operator on $X$, we have  by reflexivity that $\alpha^*= \lambda \alpha +s$ where $s$ is strictly singular.  It follows that
$\alpha=\lambda \alpha^* +s^*$, and then  $\lambda^2=1$. If $\lambda=-1$ we have an almost symplectic structure. Then by the previous proposition $X$ or its hyperplanes admit symplectic structure,  and if $\lambda=1$ we have a Hermitian  structure.

In the complex case  we may assume an isomorphism $\alpha$ of $X$ onto $\overline{X}^*$, then $\overline{\alpha}^*=\lambda \alpha +s$ and we get that $|\lambda|=1$.
If $\lambda=e^{i\theta}$  then by taking  $\mu=e^{i(\theta-\pi)/2}$ we have
$$\overline{(\mu \alpha)}^*=\overline{\mu} \overline{\alpha}^*=\overline{\mu}\lambda \alpha +s_1=-\mu \alpha +s_1,$$
for a strictly singular operator $s_1$. Hence we obtain an almost  symplectic structure on $X$.  By taking  instead $\mu= e^{i\theta/2}$ we have an Hermitian structure on $X$.
\end{proof}

When such spaces admit a symplectic structure, then it cannot be  trivial, since it would rely on writing
$X=Y \oplus Y^*$ and a space with the $\lambda I + S$- property cannot have nontrivial complemented subspaces.

\begin{prop} \label{perturbation} Let $X$ be a reflexive real Banach space. Let $\alpha: X\to X^*$ be a symplectic isomorphism and $s:X\to X^*$ be a strictly singular operator. If $\alpha +s$ is also symplectic then $\alpha $ and $\alpha +s$ are equivalent.
\end{prop}
\begin{proof}
Recall that  $\widehat{X}$ denotes the usual complexification of a real Banach space $X$. If  $T: X\to Y$ is an operator  then  $\widehat{T}$ denotes the respective induced operator from $\widehat{X}$ to $\widehat{Y}$. Let us consider the spectrum of a real operator $\tau$ as the spectrum of its complexification, and denote it by $\sigma(\widehat {\tau})$. Now set  $S= \alpha^{-1}s$, then $ \widehat{S} = \widehat{\alpha}^{-1}\widehat{s}$ is also strictly singular.  Consider $\Gamma$ a rectangular, rectifiable,  conjugation-invariant, simple closed curve, contained in the open unit disk, and such that $\Gamma \cap \sigma (\widehat S)=\emptyset$. Denote by $U$ the bounded open connected component of $\C\setminus \Gamma$, and by $V$ the unbounded open domain of  $\C\setminus \Gamma$.  Let $\widehat P$ be the spectral projection of $\widehat S$ associated to $\sigma (\widehat S) \cap U$. By a general argument (see \cite{FGEven}), $\widehat P$ is induced by a real operator $P$ on $X$.  Let also $\widehat Q$ be the spectral projection associated to $ \sigma (\widehat S)\cap V$.  Then $\widehat S= \widehat S\widehat P + \widehat S \widehat Q$.

The spectral radius of $\widehat S \widehat P$ is strictly smaller than one. Then the series $\sum_{n\geq 1} a_n(\widehat S \widehat P)^n$ converges to an operator  $\widehat R$, where $\sum_{n\geq 1}a_n z^n$ converges to $2(-1+(1+z)^{1/2})$ for every $|z|<1$. Since the coefficients of the series are reals, then $\widehat R$ is induced by  a real operator $R=\sum_{n\geq 1} a_n (SP)^n$ on $X$ which is strictly singular.  It follows that
\begin{eqnarray} \label{bhaskara}
 \widehat R + \frac{1}{4}\widehat R^2= \widehat S\widehat P.
\end{eqnarray}

We also have that $(\widehat \alpha \widehat R)^*=-\widehat \alpha \widehat R$. Indeed,  $(\widehat \alpha (\widehat S \widehat P)^n)^* = [(\widehat s \widehat P)^*\widehat \alpha^{-1*}]^n \widehat \alpha ^*=-(\widehat s\widehat P \widehat \alpha ^{-1})^n\alpha= -\widehat \alpha (\widehat \alpha ^{-1}(\widehat s\widehat P)^n= -\widehat \alpha (\widehat S\widehat P)^n$ for every $n\in \N$,  where we used that $\widehat s^*=-\widehat s$ and hence  $(\widehat s\widehat P)^*=-\widehat s\widehat P$ with the canonical identification. Therefore

\begin{align}  \label{partialsim}
\left (\widehat P+\frac{1}{2} \widehat R \right )^* \widehat \alpha \left (\widehat P+\frac{1}{2} \widehat R \right )  &=  \left ( (I+\frac{1}{2} \widehat R) \widehat P \right )^* \widehat \alpha  \left (I+\frac{1}{2} \widehat R \right ) \widehat P\\ \nonumber
&= \widehat P^*(I + \frac{1}{2}\widehat R^*)(\widehat \alpha  + \frac{1}{2}\widehat \alpha \widehat R)\widehat P\\ \nonumber
&= \widehat P ^* \left (  \widehat \alpha  + \frac{1}{2} \widehat \alpha \widehat R + \frac{1}{2}\widehat R^*\widehat \alpha + \frac{1}{4}\widehat R^*\widehat \alpha \widehat R \right ) \widehat P\\ \nonumber
&=  \widehat P^* (\widehat \alpha +  \widehat \alpha \widehat R +  \frac{1}{4}\widehat \alpha \widehat R^2)\widehat P\\   \nonumber
&=  \widehat P^*( \widehat \alpha + \widehat s) \widehat P.
\end{align}

Let $X_0=PX$ and consider $T_1 :=I +  \frac{1}{2} R$.  Since $\beta:= T_1^*\alpha T_1 = \alpha + \mathfrak s$, where $\mathfrak s$ is strictly singular and $\beta^*=-\beta$ it follows from the proof of Proposition \ref{almostsym} that for any closed subspace $Z\subseteq X$  such that $X=Z \oplus \ker \beta$  the restriction $ \beta_{|Z} : Z \to (\ker \beta)^\perp$ is symplectic.  Now observe that $T_1$ is Fredholm of index 0 and $\ker T_1\subseteq \ker \beta$. Then  we can write $X= X_1 \oplus \ker \beta$ where the restriction of $T_1$ to $X_1$ is an isomorphism onto its image.   We may assume  by Lemma \ref{subsym} that $X_1\subseteq X_0$ and  that  $\gamma= \beta_{|X_1} : X_1 \to (\ker \beta)^\perp$ is symplectic.

Let us denote by $\Omega$ and $\omega$ the symplectic forms associated to $\alpha +s$ and $\alpha$, respectively.
Equation (\ref{partialsim}) implies  that $P ^* (T_1^* \alpha T_1)P= P^* (\alpha + s) P$. Then for every $x,y\in X_1$  we have $\Omega(x,y)=  P^*(\alpha+s)P(x)(y)= \beta(x)(y)$ and then
 $X_1$ is  a symplectic subspace of $(X, \Omega)$.  Analogously, $T_1X_1$ is a symplectic subspace of $(X, \omega)$. Hence by Corollary \ref{symio}  we can write  $X=X_1 \oplus X_1^{\Omega}$ and $X= T_1X_1\oplus (T_1X_1)^{\omega}$  where $X_1^{\Omega}$  and  $(T_1X_1)^{\omega}$ are finite dimensional symplectic subspaces with the same dimension. Then there exists an isomorphism $T_2 \,: \, X_1^{\Omega} \to  (T_1X_1)^{\omega}$ such that $\omega(T_2x, T_2y)= \Omega (x,y)$ for all $x,y\in X_1^{\Omega}$.   Hence the isomorphism $T: X_1 \oplus X_1^{\Omega} \to  T_1X_1\oplus (T_1X_1)^{\omega}$ represented by the matrix $\begin{pmatrix}
  T_1 & 0\\
  0 & T_2
\end{pmatrix}$   satisfies $T^* \alpha T= \alpha+ s$.\end{proof}

With the same idea we can prove that if $\alpha$ and $\beta$ are two symplectic maps on $X$ close enough  (in the sense of the norm of operators) then $\alpha \sim \beta$.

We obtain now a result analogous to \cite[Prop. 8]{FGEven} for symplectic structures on Hilbert spaces. Recall that a bilinear map $T:X \times X \to \R$ on a Banach space is said to be compact if its associated operator $L_T: X \to X^*$ is compact. First we need the following lemma.

\begin{lemma}\label{symspectrum} The spectrum of  a  symplectic map $\alpha$ on a real Hilbert space $\mathcal H$ has only imaginary values.
\end{lemma}
\begin{proof} Let $(\widehat{ \mathcal H}, \pin{}{}_\C)$ be the complexification of $(\mathcal H, \pin{}{})$.  It follows that $\widehat{\alpha}^{*}= -\widehat{\alpha}$.
Let $\lambda \in \partial \sigma(\widehat {\alpha})$. By \cite[Corollary 3.1]{Maurey} there is a normalized sequence $(x_n)_n$ such that $\lim_{n\to \infty} (\widehat{\alpha}(x_n)-\lambda x_n )=0$.  Since $( \widehat{\alpha} - \lambda I)^{*}=  -\widehat{\alpha} -\overline{\lambda}I $ commutes with $ \widehat{\alpha}- \lambda I$ then $\| ( \widehat{\alpha} - \lambda I)^{*}(x) \| =\| ( \widehat{\alpha}- \lambda I)x\|$ for all $x\in \widehat{\mathcal H}$. Therefore $\lim_{n\to \infty} -(\lambda + \overline{\lambda})x_n = \lim_{n\to \infty} ( -\widehat{\alpha}(x_n)-\overline{\lambda} x_n )=0$ and thus $\Re(\lambda)=0$.\end{proof}

\begin{prop} \label{pertsymplectic} Let $\mathcal H$ be a real infinite-dimensional Hilbert space and $\omega$ be a symplectic structure on $\mathcal H$. Then  there do not exist bounded bilinear maps $\Omega$ and $\kappa$ on $X = \R \oplus \mathcal H$ with $\kappa$ compact, $\Omega$ an extension of $\omega$ and such that $\Omega+\kappa$ is symplectic on $X$.
\end{prop}

\begin{proof} Assume that  $\omega(x,y)=\pin{Jx}{y}$ for an equivalent  inner product $\pin{ \,}{ \,} $  on $\mathcal H$ where $J\in \mathfrak  L(\mathcal H)$ is a complex structure. $X$  equipped with the inner product determined by $\pin{(a,x)}{(b,y)}_X= ab+ \pin{x}{y}$. Let  $A$ be  the operator  on $X$ defined by the matrix
$\begin{pmatrix}
  1 & 0\\
  0 & J
\end{pmatrix}$.
Suppose that there exists a compact operator $K$ on $X$ such that $A+K$ is symplectic on  $X$ and consider the function from $[0,1]$ into $\mathfrak L(\widehat X)$ given by $T_\mu = \widehat A + \mu \widehat K$.
The operator $\widehat A$ is defined by the  matrix
$\begin{pmatrix}
  1 & 0\\
  0 & \widehat J
\end{pmatrix}$, therefore
has only one real eigenvalue, with associated spectral projection of dimension 1.
On the other hand, by Lemma \ref{symspectrum} we have that $T_1$ does not have real eigenvalues. In order to get a contradiction, the rest of the proof is then identical to  \cite[Prop. 8]{FGEven}.
\end{proof}
The spectral arguments of the proof of \cite[Prop. 8]{FGEven} do not apply when we considering symplectic maps on Banach spaces. We do not know  whether Proposition \ref{pertsymplectic} holds for general Banach spaces.\medskip

{\em Acknowledgements:} The second author thanks V. Ferenczi for his helpful comments and suggestions which led the authors to consider almost symplectic forms.

\end{document}